\documentclass{amsart}
\usepackage{graphicx,amssymb,epsfig}
\usepackage[all]{xy}
\newtheorem{theorem}{Theorem}[section]
\newtheorem{proposition}[theorem]{Proposition}
\newtheorem{lemma}[theorem]{Lemma}
\newtheorem{corollary}[theorem]{Corollary}
\newtheorem{remark}[theorem]{Remark}
\newtheorem{example}[theorem]{Example}

\newtheorem{definition}[theorem]{Definition}

\newcommand{\bth}{\begin{theorem}}
\newcommand{\bpr}{\begin{proposition}}
\newcommand{\epr}{\end{proposition}}
\newcommand{\bco}{\begin{corollary}}
\newcommand{\eco}{\end{corollary}}
\newcommand{\ble}{\begin{lemma}}
\newcommand{\ele}{\end{lemma}}
\newcommand{\bre}{\begin{remark}\rm}
\newcommand{\ere}{\end{remark}}
\newcommand{\bex}{\begin{example}\rm}
\newcommand{\eex}{\end{example}}
\newcommand{\bde}{\begin{definition}\rm}
\newcommand{\ede}{\end{definition}}

\def\la#1{\hbox to #1pc{\leftarrowfill}}
\def\ra#1{\hbox to #1pc{\rightarrowfill}}
\def\fract#1#2{\raise3pt\hbox{$ #1 \atop #2 $}}

\def\lrar{{\ra 2}}
\def\tensor{\otimes}
\def\bbh{{\mathbb H}}
\def\bbz{{\mathbb Z}}

\def\bbp{{\mathbb P}}
\def\bbr{{\mathbb R}}

\def\bbc{{\mathbb C}}
\def\bbr{{\mathbb R}}
\def\bbq{{\mathbb Q}}

\def\gr{Gr(n,m)}
\def\grass{Gr(2,2)}
\newcommand{\Fl}{Fl(1,2,4)}
\newcommand{\Hl}{Hol_1G(2,2)}

\def\rat#1{\hbox{Rat}_{#1}}
\def\bhol#1{\widetilde{\hbox{Hol}}_{#1}}
\def\brat#1{\widetilde{\hbox{Rat}}_{#1}}
\def\hol#1{\hbox{Hol}_{#1}}

\begin{document}

\title{The Space of Linear Maps into a Grassmann Manifold}
\author{S. Kallel, P. Salvatore and W. Ben Hammouda}
\address{Universit\'e des Sciences et Technologies de Lille\\
 Laboratoire Painlev\'e, U.F.R de Math\'ematiques\\
59655 Villeneuve d'Ascq, France}
\email{sadok.kallel@math.univ-lille1.fr}
\address{Dipartimento di Matematica \\ Universit\'a di Roma Tor Vergata \\
Via della Ricerca Scientifica \\ 00133 Roma, Italy }
\email{salvator@mat.uniroma2.it}

\maketitle

\begin{abstract} We show that the space of all holomorphic maps of
  degree one from the Riemann sphere into a Grassmann manifold is a
  sphere bundle over a flag manifold. Using the notions of ``kernel" and
  ``span" of a map, we
  completely identify the space of unparameterized maps as well.
  The illustrative case of maps into the quadric Grassmann manifold is
  discussed in details and the homology of the corresponding spaces computed.
\end{abstract}

\section{Introduction}

Given a complex projective variety $M$, one can consider the space
$\hol{}(M)$ of all holomorphic maps from the Riemann sphere $\bbp^1$
into it. This can be given the structure of a quasiprojective variety \cite{john}
and thus has the homotopy type of a finite CW complex.  If
we fix basepoints in all of $\bbp^1$ and $M$, then we denote by
$\rat{}(M)$ the subspace of basepoint preserving maps. The
homeomorphism type of this subspace doesn't depend on the choice of
basepoint in $\bbp^1$ and if $M$ is homogeneous, it is independent of
the choice of basepoint in $M$ as well.

The space $\hol{}(M)$ is not in general connected.  Fixing $\alpha\in
H_2(M)$, one defines $\hol{\alpha}(M)$ to be the subspace of all
morphisms $f: \bbp^1\lrar M$ with $f_*([\bbp^1])=\alpha$.  In
\cite{perrin} and for $M=G/P$, $P$ a parabolic subgroup,
sufficient conditions on $\alpha$ are given so that $\hol{\alpha}(G/P)$ is
irreducible and smooth.

When $M=\gr$ is the Grassmann manifold of $n$-dimensional complex
linear planes in $\bbc^{n+m}$, we say that a map $f: \bbp^1\lrar\gr$
has degree $d\geq 0$ if its effect on the second homology group is
multiplication by $d$.  It turns out that two maps in $\hol{}(\gr )$
are in the same connected component if and only if they have the same
degree. We write $\hol{d}(\gr )$ the component of degree $d$ maps. Of
course $\hol{0}(\gr )=\gr$ are the constant maps, while
$\hol{1}(\gr )$ are the linear maps which we study in details in this
paper. The following is a complete description of this space.

\bth\label{hol1} Let $Fl_{(1,n)}(\bbc^{n+m})= U(n+m)/U(1)\times U(n-1)\times
U(m)$ be the flag manifold of all $(1,n)$ flags in $\bbc^{n+m}$ with
projections
\begin{eqnarray*}
  p_1&:& Fl_{(1,n)}(\bbc^{n+m})\lrar \bbp^{n+m-1}\\
  p_2&:& Fl_{(1,n)}(\bbc^{n+m})\lrar \gr
\end{eqnarray*}
Then $\hol{1}(\gr )$ has as deformation retract the total space of the
sphere bundle over $Fl_{(1,n)}(\bbc^{n+m})$ associated to the rank $m$
complex vector bundle
$$p_1^*(H )\tensor p_2^*(Q)$$ where $H={\mathcal O}(1)$ is the
line bundle dual to the tautological line bundle over $\bbp^{n+m-1}$,
and $Q$ the anti-tautological bundle over $\gr$.
\end{theorem}

Some descriptions of $\hol{}(\gr )$ as an algebraic variety can be
found in \cite{john, david, kirwan} but no homotopy type has been
computed.  Note that the space $\rat{}(\gr)$ plays a significant role
in the theory of multilinear control systems. It corresponds to the
moduli space of $n$-input, $m$-output time invariant systems that are
controllable and observable \cite{hm} (and references therein).  A
complete and elegant computation of the homology of $\rat{d}(\gr )$
for all positive $d$ was obtained by Mann and Milgram
\cite{mm1,mm2}. The full space $\hol{d}(\gr )$ corresponds on the
other hand to ``singular'' or ``generalized'' state space systems
\cite{dmh}.  A more detailed study of $\hol{d}(\gr)$ for $d>1$ will be
pursued elsewhere.

Observe that $\hol{1}(\bbp^1)=PGL_2(\bbc )\simeq PU(2)$ acts
freely on $\hol{1}(\gr )$ by precomposition of maps. Similarly in the based
case, there is an action of $\rat{1}(\bbp^1)=Aff(\bbc )$ on
$\rat{1}(\gr )$. Here $Aff(\bbc )=\{z\mapsto az+b,
a\in\bbc^*,b\in\bbc\}$ is a semi-direct product and $Aff(\bbc)\simeq
S^1$.  Define $\brat{1}(\gr ):=\rat{1}(\gr )/Aff(\bbc )$ to be the
space of \textit{unparameterized} maps.  Similarly define
$\bhol{1}(\gr ):=\hol{1}(\gr )/PGL_2(\bbc )$.  Our next result
determines completely the homeomorphism type of these spaces and makes
the cute observation that an element in $\bhol{1}(\gr )$ is completely
determined by its ``span" and its ``kernel".  Following (\cite{gp},
p. 526), we define the ``kernel" and ``span" of a holomorphic map
$f:\bbp^1\lrar\gr$ to be respectively the intersection and span of all
vector subspaces $f(p)\subset\bbc^{n+m}$, $p\in\bbp^1$ (see \S\ref{unparameterized}).

\bth\label{main3} By sending a holomorphic map $f$ to the flag
$(\hbox{ker}(f)\subset \hbox{span}(f)\subset\bbc^{n+m})$, we construct
homeomorphisms
\begin{eqnarray*}
\bhol{1}(\gr )&\fract{\cong}{\lrar}& Fl_{(n-1,n+1)}(\bbc^{n+m})\\
\brat{1}(\gr )&\fract{\cong}{\lrar}& \bbp^{n-1}\times\bbp^{m-1}
\end{eqnarray*}
where $Fl_{(n-1,n+1)}(\bbc^{n+m})$ is the variety of all $(n-1,n+1)$-flags in
$\bbc^{n+m}$ of complex dimension $nm+n+m-3$.
\end{theorem}

A decent part of this paper studies the special case of holomorphic maps
into the quadric grassmann manifold $\grass$. This is the variety
realized via the Pl\"ucker embedding as a quadric hypersurface in
$\bbp^5$.  It has the homotopy and homology groups of
$\bbp^2\times S^4$ (cf. \S\ref{quadric}, \cite{lai}).  We shall prove
that

\bth\label{grasscase}\ \\
(i) There is a homotopy equivalence
$\rat{1}(\grass)\simeq S^2\times S^3$  (Corollary \ref{trivial})\\
(ii) There is a (non-multiplicative) homotopy equivalence (Proposition \ref{loopgrass})
$$\Omega (\grass )\simeq S^1\times S^3\times \Omega S^5\times\Omega S^7$$
(iii) The cohomology groups $\tilde H^*(\hol{1}(\grass );\bbz )$ are given by
(Theorem \ref{homcalc})
$$
\begin{tabular}{|c|c|c|c|c|c|c|c|c|c|c|c|c|c|c|c|c|c|c|c|}
  \hline
  $i$ & 2
  & 4 & 6 & 7 & 8 & 9 & 11 & 13 \\
  $H^i$& $\bbz\oplus\bbz$ & $\bbz\oplus\bbz$ & $\bbz_4\oplus\bbz$ & $\bbz$ & $\bbz_4$ & $\bbz\oplus\bbz$ & $\bbz\oplus\bbz$ & $\bbz$ \\ \hline
\end{tabular}$$
\end{theorem}

Part (iii) follows from a careful analysis of the Gysin sequence
associated to the bundle given in Theorem \ref{hol1} and from the
classical description of P. Baum of the cohomology of flag varieties.
This allows for effective homology computations and for a complete
description of the differentials in the evaluation fibration
$\rat{1}(\grass )\lrar\hol{1}(\grass)\fract{ev}{\lrar}\grass$ which
would have been otherwise quite difficult to obtain \cite{walid}.

\bre Write $\iota$ the inclusion of based holomorphic maps into all
based continuous maps
$$\iota : \rat{1}(\gr )\hookrightarrow\Omega^2_1(\gr )\ \ \ ,\ \ \ 1\leq n\leq m$$
The induced map in homology is determined in \cite{mm2}. This is a
non-trivial calculation and aspects of this are discussed in \S3. The
important claim made in \cite{mm2} is that the homology of
$\rat{k}(\gr )$ is generated via homology operations by the homology
of $\rat{1}(\gr )$. This point we shall return to in the future. \ere

\section{The Based Linear Maps}\label{intro}

The Grassmann manifold $\gr$ is the homogeneous space
$U(n+m)/U(n)\times U(m)$. It is a smooth complex variety of dimension
$nm$. A system of charts for $\gr$ is given as follows.  Choose a
decomposition of $\bbc^{n+m}$ into a direct sum $\bbc^n\oplus\bbc^m$,
and let $L:\bbc^n\lrar\bbc^m$ be any linear operator. Then its graph
is an $n$-dimensional subspace of $\bbc^{n+m}$. The set
$Hom(\bbc^n,\bbc^m)$ of all subspaces obtained in this way is an open
dense subspace $V$ of $\gr$. Thus any point of $V$ can be represented
by an $m\times n$ complex matrix and this is a cell of dimension $nm$.
By choosing various coordinate decompositions of $\bbc^{n+m}$, we can
cover $\gr$ with affine charts this way.

Identify $\gr$ with all $n$-planes $W$ in some decomposition $U\oplus
Y$ where $U\cong\bbc^n$ and $Y\cong\bbc^m$. Then the maximal Schubert
cell in $\gr$ is identified with $Hom(U,Y)$ and its complement is a
divisor (a Schubert hypersurface)
$$\Delta = \{W\ |\ \dim (W\cap Y)\geq 1\}$$
A based holomorphic map $f: \bbp^1=\bbc\cup\infty\lrar\gr$ can then be
thought of as a \textit{meromorphic} map into this Schubert cell
$\bbc\lrar Hom(\bbc^n,\bbc^m)$, sending $\infty$ to the trivial
matrix. The degree of such a map is computed as the intersection
number with the Schubert hypersurface; $\deg f = (f(\bbp^1):W)$.  Such
a map can be written as  a \textit{rational matrix}\footnote{Known
  in control theory as the ``transfer matrix''}
\begin{equation}\label{rational}
f(z) = \sum {A_i\over (z-z_i)^{k_i}}
\end{equation}
with $A_i$ an $m\times n$ matrix. The degree of the map in (\ref{rational})
depends on the multiplicity of the poles $z_i$ and the ranks of the
$A_i$ \cite{dh}.

\bre\label{description} There is another useful expression
for $f$ given in matrix form with entries in $\bbc [z]$.
Indeed a holomorphic map $f:\bbp^1\lrar\gr$ sending $\infty$ to
$E_n$; the plane spanned by the
\textit{first} $n$-coordinates vectors, can be represented in the form
$$z\longmapsto \hbox{span of row vectors of}\ [D:N]$$
where $[D:N]$ is an
$n\times (n+m)$ matrix, uniquely defined up to $GL_n(\bbc [z] )$-multiplication
 on the left,  with $D\in M_{n\times n}(\bbc [z])$, $N\in M_{n\times m}(\bbc [z])$
 and the degree of the determinant of $D$ is maximal among all $n\times n$ minor determinants.
 Out of this representation we can recover the transfer function described in
 (\ref{rational}) according to the formula $T(z) =
D^{-1}N$. This correspondence is discussed in various books in linear control theory (see also
\cite{walid, mm1}). We will use both representations (as a matrix form or a transfer
function) in \S\ref{unparameterized}.
\ere

Based on the transfer map description of based holomorphic maps
(\ref{rational}), we give a shorter proof of the following observation due to Mann and Milgram.

\bpr\label{rat1} $\rat{1}(Gr(n,m) )$ is homeomorphic to
$$\bbc\times (\bbc^{m}-\{0\})\times_{\bbc^*} (\bbc^n-\{0\})$$ with
$\bbc^*$ acting diagonally; i.e. $a ({\bf v}, {\bf w}) = (a^{-1}{\bf
  v}, a{\bf w})$.  This is up to homotopy the space of rank one
matrices of size $m\times n$.  \epr

\begin{proof} Elements in $\rat{1}(Gr(n,m))$ are maps of the form
  $z\mapsto {A\over z-z_0}$, where $A$ is a rank one $m\times n$
  matrix.  Any such matrix can be written as a product ${\bf v}.{\bf
    w}^T$ where ${\bf w}^T$ is a non-zero row vector of $\bbc^n$ and
  ${\bf v}$ is a non-zero column vector of $\bbc^m$.  Note that
  $(a^{-1}{\bf v},a{\bf w})$ with $a\in\bbc^*$ gives the same matrix
  $a^{-1}{\bf v}.{a}{\bf w}^T = {\bf v}.{\bf w}^T$. It is then
  immediate to see that the assignment
\begin{eqnarray*}
\rat{1}(Gr(n,m) )&\fract{\cong}{\lrar}&\bbc\times
\left((\bbc^{m}-\{0\})\times_{\bbc^*} (\bbc^n-\{0\}) \right)\\
{A\over z-z_0}&\longmapsto& (z_0, {\bf v}\times_{\bbc^*} {\bf w})
\end{eqnarray*}
is a homeomorphism, where $A = {\bf v}.{\bf w}^T$ is a rank one
$m\times n$ matrix.
\end{proof}

Proposition (\ref{rat1}) makes $\rat{1}(\gr )$ into a bundle over
$\bbp^{n-1}$ with fiber $\bbc^m-\{0\}\simeq S^{2m-1}$. We can identify
this fiber in a different way. Consider the inclusion $\iota_m:
\rat{1}(\bbp^m)\hookrightarrow\rat{1}(\gr )$ which is induced from the
inclusion $\bbp^m=Gr(1,m)\hookrightarrow\gr$. According to the above
Proposition, $\rat{1}(\bbp^m)\simeq S^{2m-1}$ and in fact
$\iota_m$ is up to homotopy the inclusion of a fiber.

\ble\label{fibration} There is a bundle
$$\rat{1}(\bbp^m)\fract{\iota_m}{\lrar}\rat{1}(\gr )\lrar\bbp^{n-1}$$
\ele

\begin{proof} Using the affine coordinate description of $\gr$
  discussed above, the inclusion $\bbp^m\hookrightarrow\gr$ is
  described at the level of a chart by the map
  $Hom(\bbc^1,\bbc^m)\lrar Hom(\bbc^{n},\bbc^m)$ which to $f:z\mapsto
  f(z)$ associates $(z_1,\ldots, z_n)\mapsto f(z_1)$. At the level of
  matrices, if $f\in \rat{1}(\bbp^m)$ is given by $z\mapsto {1\over
    z-a}[a_1,\ldots, a_m]^T$, then
$$\iota_m (f) :z\longmapsto  {1\over z-a}\begin{pmatrix}a_1&0&\cdots&0\\
  \vdots\\ a_m&0&\cdots&0\end{pmatrix}$$ But if we identify
$\rat{1}(\gr )$ with $\bbc\times
(\bbc^{m}-\{0\})\times_{\bbc^*}(\bbc^n-\{0\})$ as in Proposition
\ref{rat1}, then the map $\iota_m$ takes the form
$$a \times (a_1,\ldots, a_m)\longmapsto a \times
(1,0,\ldots, 0)\times_{\bbc^*}(a_1,\ldots, a_m)$$
which is indeed the inclusion of the fiber over $[1:0:\cdots :0]$.
\end{proof}

Next is a useful corollary to Proposition \ref{rat1}.

\bco\label{trivial}
For $n=2$ and $m$ even, $\rat{1}(Gr(2,m))\simeq S^2\times S^{2m-1}$.
\eco

\begin{proof}
  Identify $\bbc^2$ with $\bbh$ the quaternions and let it act on
  $\bbc^{m}= (\bbc^2)^{m\over 2}$ via quaternionic multiplication on
  each factor. Since $\bbh$ is a division algebra, this action
  restricts to an action of $\bbc^2-\{0\}$ onto $\bbc^m-\{0\}$.  If we
  write $\rat{1}(Gr(2,m))$ as $(\bbc^2-\{0\})\times_{\bbc^*}
  (\bbc^m-\{0\})$ as in Proposition \ref{rat1}, then this quotient is
  diffeomorphic to $\bbp^1\times (\bbc^m-\{0\})$ via the map sending
  the class $[v,x]\longmapsto ([v],v.x)$, with $v.x$ meaning $v$
  acting on $x$ via quaternionic multiplication.
\end{proof}

Corollary \ref{trivial} can also be derived from a general description
 of $\rat{1}(Gr(2,m))$ as an attachment space.
Consider the the sphere fibration obtained from Lemma \ref{fibration};
\begin{equation}\label{rat1g22}
  S^{2m-1}\lrar\rat{1}(Gr(2,m) )\lrar S^2 \ \ \ ,\ \  \ m\geq 2
\end{equation}
This bundle has a section since the euler class of its associated
bundle lies in the trivial group $H^{2m}(S^2)$.  We quickly review
some classical but structural results on sphere bundles over spheres
due to James-Whitehead and Sasao (see \cite{sasao}). We make use of a
bit of notation: (i) by an unlabeled map $S^{k-1}\lrar\Omega S^k$ we
mean the map adjoint to the identity of $S^k$, and (ii) a fibration
$E\lrar B$ with fiber $F$ has a holonomy map $h:\Omega B\lrar Aut(F)$;
with $Aut(F)$ being the self-homotopy equivalences of $F$ and which
classifies the fibration up to fiber homotopy type (Stasheff).  Let
$E$ be the total space of an $S^k$-bundle over another sphere $S^n$.
Then the composite
\begin{equation}\label{first}
f: S^{n-1}\times S^k\lrar\Omega S^n\times S^k\fract{h}{\lrar} S^k
\end{equation}
determines the bundle completely up of fiber homotopy equivalence,
where $h$ is the holonomy again. When the bundle has a section, this
map factors through the half-smash $S^{n-1}_+\wedge S^k\simeq S^k\vee
S^{n-1}\wedge S^k \lrar \Omega S^n_+\wedge S^k$ and hence defines a
map of the same name
\begin{equation}\label{themapf}
f: S^{n-1}\wedge S^k {\lrar} S^k
\end{equation}
which in turn defines by adjunction an element $[f]\in
\pi_{n-1}(Aut(S^k))=\pi_{n-1}(\Omega^k_{\pm 1}S^k)
\cong\pi_{n+k-1}(S^k)$. If $\iota_k: S^k\hookrightarrow S^k\vee S^n$
is the inclusion of the first factor, we will write $[\iota_k\circ f]$
the homotopy class in $\pi_{n+k-1}(S^k\vee S^n)$.

\bpr\label{sphereoversphere} If $E$ is an $S^k$-bundle over $S^n$ with
a section, then up to homotopy
$$E = (S^k\vee S^n)\bigcup{}_{[\iota_k,\iota_{n}]+ [\iota_{k}\circ f]}e^{n+k}$$
\epr

\begin{proof}
Write as in (\cite{stasheff}, proposition 2) the total space $E$ in the form
$$E=  D^n\times S^k\cup S^k/\ (x,y)\sim  f(x,y)\ ,\ \ (x,y)\in S^{n-1}\times S^k$$
We can then think of $E$ as the CW complex $S^k\cup_{\alpha}
e^n\cup_{\beta} e^{k+n}$ with $\alpha (x) = f(x,*)$ and $f$ as in
(\ref{first}), $\alpha'$ is the characteristic map for $e^n$;
i.e. $\alpha = \partial\alpha'$, and
$$\beta : S^{k+n-1}=D^n\times S^{k-1}\cup S^{n-1}\times D^k\lrar
D^n\times *\cup S^{n-1}\times S^k\fract{\alpha'\cup f}{\ra 3}
e^n\cup_\alpha S^k$$
If the bundle has a section, then $\alpha$ is null-homotopic and
$\beta$ becomes
\begin{equation}\label{secondmap}
D^n\times S^{k-1}\cup S^{n-1}\times D^k\lrar
D^n\times *\cup S^{n-1}\times S^k
\fract{\phi\times *\cup *\times f}{\ra 4} S^n\times *\cup *\times S^k
=S^n\vee S^k
\end{equation}
where $\phi$ is collapsing the boundary of the disk. As we pointed out and in
the presence of a section, the map $f$ factors through $S^{n-1}\wedge
S^k=S^{n+k-1}\lrar S^k$.  By construction the homotopy class of
(\ref{secondmap}) is the element $[\iota_k,\iota_n] + i_{k_*}([f])\in
\pi_{n+k-1}(S^n\vee S^k)$.
\end{proof}

We use this Proposition to completely determine the homotopy type of
$\rat{1}(Gr(2,m))$ and recover in particular Corollary \ref{trivial}.

\bco $\rat{1}(Gr(2,m))$ is up to homotopy the CW complex
$$(S^2\vee S^{2m-1})\bigcup{}_{[\iota_2,\iota_{2m-1}]+ m [\iota_{2m-1}]\circ \eta}e^{2m+1}$$
where $\eta\in \pi_{2m}(S^{2m-1})\cong\bbz_2$ is the (Hopf) generator.
\eco

\begin{proof}
  $\rat{1}(Gr(2,m))$ is up to homotopy the Borel construction
  $S^3\times_{S^1}S^{2m-1}$. The projection onto $S^2$ has a
  section\footnote{The set of all such sections is in one to one
    correspondence with $S^1$-equivariant maps of $S^3$ into
    $S^{2m-1}$.}.  The existence of the section identifies
  $\rat{1}(Gr(2,m))$ with an adjunction space $(S^2\vee
  S^{2m-1})\bigcup_{g}e^{2m+1}$ where $g: S^{2m}\fract{\alpha}{\lrar}
  S^{2m-1}\hookrightarrow S^2\vee S^{2m-1}$.  According to
  (\ref{themapf}), the map $\alpha$ is adjoint to a map
  $S^{1}\rightarrow \Omega^{2m-1}_1S^{2m-1}$. To understand this map,
  consider the diagram
$$\xymatrix{
  &&U(m)\ar[d]\ar[rd]^{J_m}&\\
  S^1\ar[r]^\iota&\Omega S^2\ar[r]^{h\ \ \ \ }\ar[ru]^{h_1}\ar[rd]_{h_2}&L^{2m-1}_1S^{2m-1}\ar[r]^\sigma&\Omega^{2m}_1S^{2m}\\
  &&\Omega^{2m-1}_1S^{2m-1}\ar[u]\ar[ru]&\\
}$$ where $h$ is the holonomy of the bundle (\ref{rat1g22}), which
factors through $h_2$, $h_1$ is the holonomy of the bundle $m{\mathcal
  O}(-1)$ of which (\ref{rat1g22}) is the sphere bundle, $\sigma$ is
suspension and finally $J_m$ is the map given by considering $U(m)$ as
self-transformations of $\bbc^m$ then compactifying.  The map $\alpha
: S^{2m}\lrar S^{2m-1}$ is by construction the adjoint to
$h_2\circ\iota$. We claim that $\alpha$ is up to homotopy $-m\eta$
where $\eta$ is the representative of the Hopf generator in
$\pi_{2m}(S^{2m-1})$.  To see this, it is enough to show that the
composite $J_m\circ h_1\circ\iota$ is adjoint to $-m\eta$ as well.

The bundle $m{\mathcal O}(-1)$ over $S^2$ is classified by a
\textit{clutching number} in $\pi_1(U(m))=\bbz$ represented by
$h_1\circ\iota$. This clutching number is $-m\in\bbz$ since this
clutching map is the composite $S^1\rightarrow U(1)^m\rightarrow U(m)$
sending
$$\lambda\longmapsto (\lambda^{-1},\cdots , \lambda^{-1})\longmapsto
\begin{pmatrix}\lambda^{-1}&&0\\
  &\because&\\
  0&&\lambda^{-1}\end{pmatrix}$$ The composite of this map with the
determinant map $U(m)\rightarrow S^1$ sends
$\lambda\longmapsto\lambda^{-m}$ from which we deduce that the
clutching number is $-m$. So what remains to be seen is that $J_m$ sends
the generator of $\pi_1(U(m))$ to the hopf generator in
$\pi_{2m+1}(S^{2m})$. But this is a direct consequence of the diagram
$$\xymatrix{
  U(1)=S^1\ar[r]^{J_1}\ar[d]^\subset&\Omega^2_1S^2\ar[d]^\sigma\\
  U(m)\ar[r]^{J_m}&\Omega^{2m}_1S^{2m} }$$ and the fact that the
adjoint of $J_1$ can be seen to be a Hopf generator by looking at
the linking number.
\end{proof}

\bpr Write $y,z$ $H^*(\rat{1}(Gr(2,m);\bbz_2)$ the generators in
degrees $2m-1$ and $2m+1$ respectively. Then $Sq^2(y) = mz$, with
$Sq^2$ is the Steenrod squaring operation.  \epr

\begin{proof} This comes down to showing that for $m$ odd,
  $Sq^2y=z$. But this is equivalent to having the attaching map of the
  $2m+1$-cell in $\pi_{2m}(S^2 \vee S^{2-1})$ project onto the Hopf
  class in $\pi_{2m}(S^{2m-1})$.
\end{proof}

\section{The Quadric Grassmannian}

The Grassmann variety $\grass$ can be realized via the Plucker
embedding $\wp : \grass\hookrightarrow\bbp^5$ as a hypersurface of
degree $2$ given in homogeneous coordinates by
$z_0z_1+z_2z_3+z_4z_5=0$.  There is a general fact we now recall. If
$X\subset\bbp^n$ is a nonsingular hypersurface of degree $d$, then its
diffeomorphism type is determined completely by $n$ and $d$. In fact
two such hypersurfaces $X,Y$ in $\bbp^n$ are ambiantly isotopic in the
sense that there is a diffeomorphism of $\bbp^n$, isotopic to the
identity, which restricts to give a diffeomorphism of $X$ to $Y$ (see
\cite{kw}, \S4).  Consequently we can choose the hyperquadric defined
by the``Fermat" equation
$$Q_{n} := \{z_0^2 +\cdots  + z_{n+1}^2 =0\}$$
as our model hypersurface (this is simply connected if $n\geq 1$).
Since $\grass$ is also a hyperquadric in $\bbp^5$, it is diffeomorphic
to $Q_4$. This ambiant diffeomorphism being needed later is made
explicit below.

\ble\label{self} The selfmap $\gamma: \bbp^5\rightarrow\bbp^5$ which takes
$[z_0:\ldots : z_5]$ to
$$
\left[ {z_1-z_0\over 2}, i{(z_0+z_1)\over 2},
{z_3-z_2\over 2}, i{(z_2+z_3)\over 2},
{z_5-z_4\over 2}, i{(z_4+z_5)\over 2} \right]
$$
restricts to a diffeomorphism between $\grass$ and $Q_4$.
\ele

\bre An abstract homeomorphism $Q_4\cong\grass$ can be obtained as
follows.  Let $G^+(2,2n)\cong SO(2n+2)/SO(2)\times SO(2n)$ be the real
grassmann variety of oriented $2$-planes in $\bbr^{2n+2}$. Then there
is a diffeomorphism $G^+(2,2n)\cong Q_{2n}$, for all $n\geq 1$, which
sends an oriented $2$-plane spanned by orthonormal vectors $v_1,v_2$
to $\pi (v_1+iv_2)$ in $Q_{2n}$, where $\pi$ is the natural projection
$\bbc^{2n+2}\backslash\{0\}\lrar\bbp^{2n+1}$ (eg. \cite{lai}).  It remains to
identify $\grass$ with $G^+(2,4)$
the real Grassmann manifold of all
oriented two plane subspaces of $\bbr^6$, and this can be done through
the sequence of homeomorphisms:
\begin{eqnarray*}
\grass := SU(4)/S(U(2)\times U(2))
&\cong& Spin(6)/(U(1)\times SU(2)\times SU(2))\\
&\cong& Spin(6)/(SO(2)\times Spin(4))\\
&\cong& SO(6)/(SO(2)\times SO(4)) =: G^+(2,4)
\end{eqnarray*}
using the standard group isomorphisms $U(1)\cong SO(2), SU(4)\cong Spin(6)$ and
$Spin(4)\cong SU(2) \times SU(2)$.
Notice the cute result that $Q_4$, and hence
$\grass$, has the same homology and same homotopy groups as
$\bbp^2\times S^4$ according to \cite{lai}.\ere

In all cases and as a consequence of the above we can derive the
following result.

\bpr\label{loopgrass} There is a non-multiplicative splitting
$$
  \Omega\grass \simeq S^1\times S^3\times\Omega S^5\times\Omega S^7
$$
\epr

\begin{proof} By non-multiplicative we mean that there is no $H$-map
  from the the left to the righthand side which is a homotopy
  equivalence. The $H$-space structure on the righthand side is the
  obvious product of $H$-space structures. The proof of this Theorem
  is based on an observation from \cite{ps}.  Replace $\grass$ by the
  Fermat hypersurface $Q_4 = \{\sum z_i^2 = 0\}$ and consider the
  pullback of the Hopf fibration
\begin{equation}\label{pullback}
\xymatrix{\tilde Q_4\ar[d]\ar[r]&S^{11}\ar[d]\\
Q_4\ar@{^(->}[r]&\bbp^5}
\end{equation}
Write $u_i = Re(z_i)$
and $v_i=Im(z_i)$. Then
$$S^{11} = \left\{(u_1,v_1,\ldots, u_6,v_6)|\ \sum u_j^2+\sum v_j^2 = 2\right\}$$
The pullback is given by
$$\tilde Q_4 = \left\{ (u_1,v_1,\ldots, u_6,v_6)\ |\ \sum u_j^2=\sum v_j^2 = 1\ |\
  \sum u_jv_j = 0\right\}$$ and this is evidently homeomorphic to the space
  of orthonormal 2-frames in $R^6$, or equivalently to the
  unit tangent bundle $ST(S^5)$ to $S^5$. This is an $S^1$-bundle over $\grass$ and
so we therefore have the fibration
\begin{equation}\label{thefib}
\Omega ST(S^5)\lrar\Omega\grass\fract{p}{\lrar} S^1
\end{equation}
Note that the projection $p$ induces an isomorphism on fundamental
groups as can easily be checked.  The fibration (\ref{thefib}) has a
homotopy retract $S^1\rightarrow\Omega\grass$. A justification of this
fact goes as follows: suppose $X$ is a simply connected space with
$\pi_2(X)=\bbz$. Represent the adjoint of the generator in $\pi_2(X)$
by a map $\alpha : S^1\lrar\Omega X$. Since $\pi_1(\Omega X)
\cong\bbz$, then $H^1(\Omega X)\cong\bbz$ and this is represented by a
map $p: \Omega X\lrar S^1$ which is necessarily a homotopy retract to
$\alpha$ (\footnote{This also follows from the following fact in
  homotopy theory: Let $\Omega E\lrar\Omega B\fract{p}{\lrar}
  F\fract{i}{\lrar} E\fract{f}{\lrar} B$ be the sequence of fibrations
  induced from $f$. Then $i$ is homotopically trivial if and only if
  $p$ has a homotopy section.}).

In all cases, the existence of a homotopy section for the multiplicative fibration
(\ref{thefib}) means that $\Omega\grass\simeq S^1\times\Omega
ST(S^5)$. On the other hand, $ST(S^5)\lrar S^5$ also has a section and
thus
\begin{eqnarray}\label{decomp}
  \Omega\grass\simeq S^1\times\Omega ST(S^5)&\simeq&
  S^1\times\Omega S^4\times\Omega S^5\nonumber\\
  &\simeq&S^1\times S^3\times\Omega S^7\times\Omega S^5
\end{eqnarray}
Here we use as well the fact that $\Omega S^4\simeq S^3\times\Omega
S^7$.

To show that our splitting is not multiplicative, we look at the loop
homology algebra of $\grass$ and show that it differs from the loop
homology algebra of the righthand side of the decomposition.  To that
end one uses the fact that if $X$ is a simply connected space of
finite type that is \textit{formal} over $\bbz$, then its cohomology
algebra determines completely its Pontryagin loop algebra
structure. To explain what this means, a space is formal if its
cohomology algebra is quasi-isomorphic to its singular cochain algebra. This
means that there are quasi-isomorphisms of differential graded
(associative) algebras
$$(C^*(X),d)\longleftarrow (A,d)\lrar (H^*(X),0)$$
There is a very useful criterion given by \cite{haouari} for showing
when a simply connected finite type space $X$ is formal when working with
coefficients over a field or over $\bbz$ if $X$ has torsion free homology. This consists
in showing that the minimal Adams-Hilton model for $X$ has a purely quadratic
differential. We recall that the Adams-Hilton model is a differential
graded algebra of the form $TV$, a tensor algebra on $V=s^{-1}\tilde
H_*(X)$ (coefficients in a commutative ring and $s$ is the suspension
operator) with differential $d$ which decomposes as $d_2+d_3+...$,
where $d_k$ is the part that maps into length-$k$ decomposables. We
say that $TV$ has a purely quadratic differential if $d_k=0$, $k>2$.

On the other hand for any simply connected $X$, we have the
quasi-isomorphism of algebras
$$C_*(\Omega X)\simeq \underline{\Omega}(C_*(X))\simeq (BC^*(X))^\vee
$$
where $\underline{\Omega}$ is Adams cobar construction, $B$ is the
bar construction, $^\vee$ the hom-dual and
$C_*(X)$ are the singular chains. The Pontryagin ring structure only depends
on the algebra structure on $C^*(X)$, so that if $X$ is formal we can replace
$C^*(X)$ by $H^*(X)$ as algebras. In other words , for formal spaces
the ring structure of $H_*(\Omega X)$ depends only on the cohomology ring $H^*(X)$.
In particular if a formal space $X$ has the cohomology algebra of a product of spheres,
then its loop homology is the loop homology of that product of
spheres.

In the case of $ST(S^5)$, it is clear that for its minimal model, $V =
T(e_3,e_4,e_{8})$ and $d(e_{8})$ is necessarily purely quadratic (in
fact $d(e_8) = [e_3,e_4]$). It is then clear that $ST(S^5)$ is
formal. On the other hand, an easy argument using the Leray-Hirsch Theorem
 shows that $H^*(ST(S^5))\cong H^*(S^4\times S^5)$ as algebras.
This implies as we discussed that
$$H_*(\Omega ST(S^5))\cong H_*(\Omega S^4)\tensor H_*(\Omega S^5)$$
as Pontryagin rings. But $H_*(\Omega S^4)\cong T(a_3)$ is obviously
not isomorphic as algebras to $H_*(S^3\times\Omega S^7)$, because there is an exterior 3-dimensional generator only for the latter, so that the decomposition in (\ref{decomp}) cannot be one of $H$-spaces.
\end{proof}

\bre There is a bundle $U(2)\lrar V(2,2)\lrar \grass$, where $V(2,2)$
is the Stiefel manifold of orthonormal $2$-frames in $\bbc^4$.  One
can recover the decomposition  (\ref{decomp}) from this
bundle and the known homotopy equivalence $V(2,2)\simeq S^5\times
S^7$. Note that the argument of proof above shows that the projection
$\Omega\grass\lrar U(2)$ cannot have a multiplicative section \cite{walid}. \ere

\bre Many questions remain open. We suspect strongly that $\Omega ST(S^5)\simeq \Omega
S^4\times\Omega S^5$ and $\Omega\grass\simeq S^1\times\Omega
S^4\times\Omega S^5$ as \textit{$H$-spaces}.  Interestingly the induced splitting
$$\Omega^2_0\grass\simeq\Omega S^3\times\Omega^2S^5\times\Omega^2S^7$$
is \textit{not} a splitting of double loop spaces \cite{walid}. To our knowledge, the
structure of $H_*(\Omega^2\grass;\bbz_p)$ as a module over the Dyer-Lashof algebra
is not generally known. \ere

\subsection{The Homology Embedding $\iota_*$} Using the pullback
bundle (\ref{pullback}), we can try to take another view on the
geometry of the inclusion
$$\iota : \rat{1}(\grass )\hookrightarrow\Omega^2_1(\grass )$$
According to Corollary \ref{trivial}, $\rat{1}(\grass )\simeq
S^2\times S^3$ while
\begin{equation}\label{homeq}
  \Omega^2_0\grass\simeq\Omega^2ST(S^5)\simeq\Omega^2S^4\times\Omega^2S^5\simeq
  \Omega S^3\times\Omega S^7\times\Omega^2S^5
\end{equation}
Write $\alpha ,\beta$ the homology classes of degree $2,3$
respectively in $H_*(\rat{1}(\grass ))\cong H_*(S^2\times S^3)$, and
we let $a,b$ and $f$ the bottom homology generators of degree $2,3$
and $5$ in $H_*(\Omega^2_0\grass )$ according to the decomposition
(\ref{homeq}).  View $S^5\subset\bbr^6$ and let $s: S^5\lrar ST(S^5)$
be the section sending locally
$$v=(v_1,\ldots, v_6)\mapsto (v, v^\perp)\ ,\
v^\perp = (v_2,-v_1, v_4,-v_3, v_6,-v_5)$$

\bpr The following diagram commutes
$$\xymatrix{
  S^{3}\ar[r]^\simeq\ar[d]&\rat{1}(\bbp^2)\ar[r]^{\iota_2}\ar[d]&\rat{1}(\grass )\ar[d]^\iota\\
  \Omega^2S^5\ar[r]^\simeq&\Omega^2_1(\bbp^2)\ar[r]^{\iota_2}&\Omega^2_1(\grass
  )\ar[r]& \Omega^2ST(S^5) }$$ and the point is that the the bottom
composite is homotopic to the section $\Omega^2s$.  As a consequence
$i_*(\beta ) = b$ in homology.  \epr

\begin{proof}
The claim immediately follows from Lemma \ref{fibration} and
the existence of the commutative diagram
$$\xymatrix{S^1\ar[d]\ar[r]^=&S^1\ar[d]\ar[r]^=&S^1\ar[d]\\
  S^5\ar[r]^s\ar[d]&ST(S^5)\ar[r]\ar[d]^\pi&S^{11}\ar[d]\\
  \bbp^2\ar[r]^\beta\ar[rd]&Q_4\ar[r]\ar[d]^\cong&\bbp^5\\
  &\grass\ar[r]^\wp&\bbp^5\ar[u]^{\gamma'} }$$ where the top part of
the diagram is made of circle fibrations, where $\gamma'$ is a
self-diffeomorphism, $ST(S^5)$ is the sphere tangent bundle identified
with the pullback of the Hopf fibration over $\bbp^5$, $Q_4$ is the
Fermat hypersurface ambiantly diffeomorphic to $\grass$ via $\gamma$
(Lemma \ref{self}), and $\beta : [w_1,w_2,w_3]\mapsto [\bar w_1: i\bar
w_1: \bar w_2 : i\bar w_2 : \bar w_3, i\bar w_3]$. The expression of
$\beta$ comes from the fact found in the proof of Proposition
\ref{loopgrass} that if $v = (v_1,\ldots, v_6)\in S^5$, then
$$\pi (v,v^\perp) =  \left[
  v_1-iv_2 : v_2+iv_1 : v_3-iv_4 : v_4+iv_3 : v_5-iv_6 :
  v_6+iv_5\right]$$ An inspection shows that all maps in the upper
part of the diagram commute.  Remains to see the bottom part of the
diagram and what $\gamma'$ is.  The composite
$\bbp^2\hookrightarrow\grass\fract{\wp}{\lrar} \bbp^5$ is the standard
embedding $[w_1:w_2:w_3]\mapsto [w_1:0:w_2:0:w_3:0]$.  If we apply
$\gamma$ to this we obtain
$$[w_1:w_2:w_3]\longmapsto [iw_1: w_1: iw_2:w_2:iw_3:w_3]$$
This differs by an obvious self-diffeomorphism $\alpha$ of $\bbp^5$
from the top map $\beta$ explicited above. We can then set $\gamma' =
\gamma\circ\alpha$.
\end{proof}

We can also trace the effect of the map $\iota$ on $H_2$.  If we write
$\rat{1}(\grass )\simeq S^3\times_{S^1}S^3$, then as in the proof of
Lemma \ref{trivial} there is an equivalence $S^3\times_{S^1}S^3\lrar
S^2\times S^3$ sending $[a,b]\mapsto ([a], ab)$. There is then a
section $S^2\lrar S^3\times_{S^1}S^3, [a]\mapsto [a,a^{-1}]$ where
$a\in S^3$ and $[a]$ its class in $S^2=\bbp^1$. Consider the bundle
$U(2)\lrar V(2,2)\lrar\grass$ with holonomy $h : \Omega\grass\lrar
U(2)$. Since $U(2)\simeq S^1\times SU(2)$, we write
$\Omega_cU(2)\simeq\Omega SU(2)$ a component of $\Omega U(2)$.  The
following is a consequence of a calculation in (\cite{mm1}, equation
4.8).

\bpr\label{i(a)} The composite
$$
\bbp^1\fract{s}{\lrar} \rat{1}(\grass )\fract{\iota}{\lrar}\Omega^2_1(\grass )
\fract{\Omega h}{\lrar}\Omega SU(2)=\Omega S^3$$ is up to homotopy the
adjoint to the identity of $S^3$.  In homology, this means that
$\iota_*(\alpha ) = a$.
\epr

\section{Proof of Theorem \ref{hol1}}

The tautological bundle $\gamma_n$ over $\gr$ is the $n$-dimensional
complex vector bundle with total space
$$\{(P,v)\in \gr\times\bbc^{n+m}\ |\ P\in\gr, v\in P\}$$
Define ${\mathcal O}(d)$ to be the complex line bundle over $\bbp^n$
whose total space is the borel construction
$$\bbc\times_{\bbc^*} (\bbc^{n+1}-\{0\}) := \bbc\times
(\bbc^{n+1}-\{0\})/_\sim$$ which is the quotient obtained by
identifying tuples
$$(x, z_0,\ldots, z_n) \sim (\lambda^d x , \lambda z_0,\ldots, \lambda
z_n) \ \ \ \lambda\in\bbc^*, d\in\bbz
$$ Clearly the trivial line bundle is $\epsilon = {\mathcal O}(0)$. It
can be checked that
\begin{itemize}
\item $Hom_\bbc({\mathcal O}(n),\epsilon ) = {\mathcal
  O}(-n)$, and that
\item ${\mathcal O}(d_1)\tensor {\mathcal O}(d_2)\cong
{\mathcal O}(d_1+d_2)$,
more particularly ${\mathcal O}(\pm d):={\mathcal O}(\pm 1)^{\tensor d}$.
\item The tautological bundle $\gamma_1$ over $\bbp^n$ is isomorphic
  to ${\mathcal O}(-1)$. This bundle has \textit{no global} non-zero
  sections. The bundle $\gamma_1$ over $\bbr P^1$ the real Grassmann
  manifold $Gr(1,1)$ is homeomorphic to the Mobius band. The bundle
  $\gamma_1$ is called the \textit{Hopf} line bundle since the
  complement of its zero section is the bundle $\gamma_1^*$ with total
  space
$$\bbc^*\times_{\bbc^*} (\bbc^{n+1}-\{0\}) = \bbc^{n+1}-\{0\}\lrar
\bbp^n$$ and this is precisely the Hopf fibration.
\item The dual of the Hopf line bundle is the
\textit{hyperplane} bundle $H$. This is the bundle obtained by
projecting
$\bbp^{n+1} - \{[0:\cdots : 0:1]\}\lrar \bbp^n$ sending
$[z_0:\cdots z_n:z]\mapsto [z_0:\cdots :z_n]$. It can easily be seen that
$H\cong {\mathcal O}(1)$ and is the dual to the Hopf line bundle.
\end{itemize}

\bde If $\eta$ is a bundle over $\bbp^n$, we write $n\eta$ its
$n$-fold Whitney sum and write $\eta^*$ the complement of its zero
section. \ede

\begin{proof} (of Theorem \ref{hol1}) The linear action of $U(n+m)$ on
  $\bbc^{n+m}$ induces an action on $\gr$ and hence an action on
  $\hol{k}(\gr )$ by postcomposition. As in \S\ref{intro}, if we
  choose a decomposition $\bbc^{n+m}\cong U\oplus Y$ where
  $U\cong\bbc^n$ and $Y\cong\bbc^m$, then the stabilizer of this
  decomposition is a copy of $U(n)\times U(m)$ embedded in standard
  diagonal way in $U(n+m)$. The action of $U(n+m)$ on $\hol{}(\gr )$
  being transitive, we have a Borel type description
$$\hol{}(\gr ) = \rat{}(\gr )\times_{U(n)\times U(m)}U(n+m)$$
If we write a map $f\in\rat{k}(\gr )$ as a rational function
(\ref{rational}), then the action of $X\times Y\in U(n)\times U(m)$ on
the $A_i\in Hom(U,Y)$ is given by the product of matrices
$$(X,Y)\cdot A_i = Y^{-1}A_iX$$
If we identify $\rat{1}(\gr )$ with $\bbc\times
(\bbc^m-\{0\})\times_{\bbc^*} (\bbc^n-\{0\})$ as in Proposition
\ref{rat1}, the action of $U(n)\times U(m)$ translates into a diagonal
action on this product by sending
\begin{equation}\label{action}
(X,Y)\times (a, {\bf v}, {\bf w})\longmapsto (a, (Y^{-1}{\bf v},
X^T{\bf w}))
\end{equation}
with $X\in U(n), Y\in U(m), {\bf v}\in\bbc^m, {\bf w}\in\bbc^n$.

We can then rewrite up to homotopy
\begin{eqnarray}
  \hol{1}(\gr ) &=& \rat{1}(\gr )\times_{U(n)\times U(m)}U(n+m)\nonumber\\
  &\simeq&
  \left[(\bbc^n-\{0\})\times_{\bbc^*}(\bbc^m-\{0\})\right]\times_{U(n)\times
    U(m)}U(n+m)
\label{sphereb}
\end{eqnarray}
where both the actions of $U(n)$ and $U(m)$ on $(\bbc^n-\{0\})$ and
$(\bbc^m-\{0\})$ are via multiplication on the right. This is not
quite the same action as in (\ref{action}) but it yields homeomorphic
quotients. Define
\begin{eqnarray*}
X(n,m) &:=&
\left[\bbc^m\times_{\bbc^*}(\bbc^n-\{0\})\right]\times_{U(m)\times
U(n)}U(n+m)\\
&=& \bbc^m\times_{\bbc^*\times U(m)}\left[(\bbc^n-\{0\})\times_{
U(n)}U(n+m)\right]\\
&\simeq& \bbc^m\times_{U(1)\times U(m)}\left[{U(n)\over
U(n-1)}\times_{U(n)}U(n+m)\right]\\
&=& \bbc^m\times_{U(1)\times U(m)}\left[{U(n+m)\over U(n-1)}\right]
\end{eqnarray*}
where here we have replaced $\bbc^n-\{0\}$ up to equivariant
homotopy by $S^{2n-1}=U(n)/U(n-1)$ and $\bbc^*$ by $U(1)=S^1$.
The projection
\begin{equation}\label{bundle}
X(n,m)\lrar U(m+n)/U(1)\times U(n-1)\times U(m)
\end{equation}
makes $X(n,m)$ into an $m$-dimensional complex vector bundle over
$Fl_{(1,n)}(\bbc^{n+m})$ of which (\ref{sphereb}) is the sphere
bundle.  Our aim is therefore to show that (\ref{bundle}) is a bundle
isomorphic to $p_1^*(H )\tensor p_2^*(Q)$.

The action of $\lambda\in U(1)$ on $S^{2n-1}\subset\bbc^n$ is via
multiplication by $\lambda$. This means that the action of the same
$\lambda$ on $\bbc^m$ is via multiplication by $\lambda^{-1}$.  That
is, the action of $(\lambda,A)\in U(1)\times U(m)$ on $\vec
w\in\bbc^m$ is $\lambda^{-1} A({\bf w})$, and the action of $U(1)$ on
$U(n+m)/U(n-1)$ is via left multiplication if we view $U(1)$ as the
top left standard matrix subgroup of $U(n+m)$.

Note that multiplication $U(1)\times U(n-1) \lrar U(n)$ induces the
projection $p_2: Fl_{(1,n)}(\bbc^{n+m})\rightarrow \gr$ and the total
space of the pullback of $Q$ is precisely
\begin{equation}\label{bundle1}
\bbc^m\times_{U(m)}\left[U(n+m)/U(1)\times U(n-1)\right]
\end{equation}
with $U(m)$ acting via multiplication on the left.
Similarly the total space of the pullback of the
dual hyperplane bundle via the projection $p_1$ induced from
$U(n-1)\times U(m) \rightarrow U(n+m-1)$ is
\begin{equation}\label{bundle2}
\bbc\times_{U(1)}\left[U(n+m)/ U(n-1)\times U(m)\right]
\end{equation}
with $U(1)$ acting by multiplication by the inverse.  We rewrite both
bundles (\ref{bundle1}) and (\ref{bundle2}) as in
\begin{eqnarray*}
&\bbc^m\times_{U(1)\times U(m)}\left[U(n+m)/U(n-1)\right]\\
&\bbc\times_{U(1)\times U(m)}\left[U(n+m)/U(n-1)\right]
\end{eqnarray*}
So that their tensor product is the bundle
$$ (\bbc\tensor\bbc^m)\times_{U(1)\times
  U(m)}\left[U(n+m)/U(n-1)\right]
$$
with $(\lambda, A)\in U(1)\times U(m)$ acting on $z\tensor {\bf w}$ by
$\lambda^{-1}z\tensor A{\bf w}$. Upon identifying $\bbc\tensor\bbc^m$
with $\bbc^m$, we recover precisely the bundle $X_{n,m}$ and the
theorem is proved.

The proof of Theorem \ref{hol1} shows that for $n\leq m$ there is a
diagram of fibrations
$$\xymatrix{
S^{2m-1}\ar[r]^{=}\ar[d]&S^{2m-1}\ar[d]&\\
Rat_1(\gr )\ar[r]\ar[d]&Hol_1(\gr )\ar[d]\ar[r]^{ev}& \gr\ar[d]^=\\
\bbp^{n-1}\ar[r]& Fl(1,n,n+m)\ar[r]& \gr
}$$
where the left vertical fibration is same as the sphere bundle of
$m{\mathcal O}(-1)$ over $\bbp^{n-1}$.
\end{proof}

We can now derive a few interesting corollaries.

\bco\label{unitsphere} $\hol{1}(\gr )$ is of the homotopy type of a
closed oriented manifold of dimension $2m-1 + \dim Fl_{(1,2)}(\bbc^{n+m}) =
2n(m+1)+2m-3$.  \eco

\bco\cite{ks}\ $\hol{1}(\bbp^m)$ is up to homotopy the unit tangent
bundle of $\bbp^m$. \eco

\begin{proof} In this case $Fl_{(1,1)}(\bbc^{1+m})=\bbp^m$ and
  $\hol{1}$ is up to homotopy the total space of the sphere bundle of
  $Q^\vee\tensor\gamma$. Here $Q^\vee$ is the
  notation for the dual of the anti-tautological bundle $Q$. We need
  check that this tensor product is isomorphic to the tangent bundle
  of $\bbp^n$. We know that $Q\oplus\gamma \cong (n+1)\epsilon$ the
  trivial bundle of rank $(n+1)$ over $\bbp^n$ and hence by taking
  duals $Q^\vee\oplus\gamma^\vee\cong (n+1)\epsilon$. We can tensor
  both sides by $\gamma$ and get
$$(Q^\vee\tensor\gamma) \oplus \epsilon = (n+1)\gamma$$
using the fact that $\gamma^\vee\tensor\gamma = \epsilon$.
On the other hand, it is quite well-known that
$$T\bbp^n\oplus\epsilon \cong (n+1)\gamma$$ This means that $T\bbp^n$
and $Q^\vee\tensor\gamma$ are stably isomorphic. But within this range
the two bundles must be isomorphic according to \cite{husemoller},
chapter 9, Theorem 1.5.
\end{proof}

The next corollary recovers a calculation of
Cohen-Lupercio-Segal \cite{cls} who identify
$\rat{1}(BU(n))$ with $F_{1,n}\simeq\bbp^{n-1}$ the first piece in the Mitchell
filtration of $\Omega SU(n)$. Here $\rat{1}(BU(n))$ is defined as the direct limit
of the inclusions $\rat{1}(\gr )\hookrightarrow\rat{1}(Gr(n,m+1))$.

\bco\label{clsn1} There are homotopy equivalences
$\hol{1}(BU(n))\simeq Fl(1,n)(\bbc^\infty )$ and
$\rat{1}(BU(n))\simeq\bbp^{n-1}$. Moreover
$H^*(\hol{1}(BU(n))$ is a free $H^*(BU(n) )$-module
with generators $1,\zeta,\ldots, \zeta^{n-1}$.
\eco

\begin{proof} Theorem \ref{hol1} gives a fibration
$$S^{2m-1}\lrar\hol{1}(\gr )\lrar Fl_{(1,n)}(\bbc^{n+m})$$
which we can stabilize via the embeddings $Fl_{(1,n)}(\bbc^{n+m})\lrar
Fl_{(1,n)}(\bbc^{n+m+1})$ that take a $(1,n)$-flag in $\bbc^{n+m}$ and
embed it in $\bbc^{n+m+1}$. In the limit we get the diagram of
fibrations
$$\xymatrix{
   S^\infty\ar[r]\ar[d]^=&\rat{1}(BU(n))\ar[r]\ar[d]&\bbp^{n-1}\ar[d]\\
  S^\infty\ar[r]\ar[d]&\hol{1}(BU(n))\ar[r]\ar[d]& Fl_{(1,n)}(\bbc^{\infty})\ar[d]\\
  \cdot\ar[r]&BU(n)\ar[r]^\cong&Gr(n,\infty) }$$ and the first two claims follow
at once since $S^\infty$ is contractible.  The cohomological
calculation is on the other hand a general consequence of the
Leray-Hirsh theorem (see for example Switzer's book, Theorem
15.47). If $\bbp^k\fract{j}{\lrar} E\lrar B$ is a fibration with an
element $\zeta\in H^2(E)$ such that $j^*(\zeta )$ maps to the
generator of $H^2(\bbp^k)$, and thus $j^*$ is an epimorphism, then
$H^*(E)$ is a free $H^*(B)$-module with generators $1,\zeta,\ldots,
\zeta^k$.
\end{proof}

\bre The Mitchell filtration consists of compact complex subvarieties
$$F_{1,n}\subset F_{2,n}\subset\cdots F_{k,n}\subset\cdots\subset
F_{\infty,n}=\Omega_{pol}SU(n)\simeq \Omega SU(n)
$$
where $F_{k,n}$ is homeomorphic to the space of polynomial loops
$\gamma : S^1\lrar SU(n)$ such that $\gamma (1)=1$ and $\gamma (z) =
\sum_0^kA_iz^i$ where $A_i$ are $n\times n$ matrices.  In
\cite{segal2} it is stated that $F_{1,n}$ corresponds to the subspace
generated by all transformations $\lambda_V:
z\lrar \begin{pmatrix}z&0\\1&0\end{pmatrix}$ with matrix written in
terms of the decomposition $\bbc^n = V\oplus V^\perp$, with $\dim
V=1$.  This is a copy of $\bbp^{n-1}$.  Note that the inclusion
$F_{1,n}\hookrightarrow\Omega_{pol}SU(n)$ is up to homotopy the
standard map $\bbp^{n-1}\lrar\Omega SU(n)$ which generates in homology
the ring $H_*(\Omega SU(n))$.  \ere

\section{The Unparameterized Linear Maps}\label{unparameterized}

In this section we prove Theorem \ref{main3}.  We defined in the
introduction the spaces $\bhol{1}(\gr )$ and $\brat{1}(\gr )$ of
unparameterized maps.  Recall that we can associate to a morphism $f:
\bbp^1\rightarrow\gr$ its \textit{kernel}
$\hbox{ker}(f) = \bigcap_{p\in\bbp^1}f(p)$. Similarly
the \textit{span} of $f$ is the linear span of these subspaces.
Naturally ker$(f)\subset \hbox{span}(f)\subset\bbc^{n+m}$. Notice that
both \textit{ker} and \textit{span} are invariant under the action of
$PGL_2(\bbc )\simeq PU(2)$ on $\hol{1}(\gr )$.

Consider the based case first. Recall any such map is of the form
$z\mapsto f(z) = {A\over z-a}$ with $a\in\bbc$ and $A$ is an $n\times m$-matrix
of rank one which we write $A = {\bf v}.{\bf w}^T$ for some non-zero ${\bf w}\in\bbc^m$
and ${\bf v}\in\bbc^n$. As in Proposition \ref{rat1}, we have the
identification
\begin{equation}\label{identification}
\rat{1}(\gr)\cong \bbc\times (\bbc^{n}-\{0\})\times_{\bbc^*} (\bbc^m-\{0\})
\end{equation}
Write $(\alpha,\beta)\in\bbc^*\times\bbc$ the element of $Aff(\bbc )$
acting on $\bbp^1=\bbc\cup\{\infty\}$ by fixing $\infty$ and sending
$z\mapsto \alpha z+\beta$. The action of $Aff(\bbc )$ on $\rat{1}(\gr
)$ by precomposition of maps translates under the above
identifications to the action
$$(\alpha, \beta)\times (a, [{\bf v}, {\bf w}])\longmapsto
\left({a-\beta\over\alpha}, [{1\over\alpha}{\bf v}, {\bf w}]\right)
$$
where again $(\alpha, \beta)$ is viewed as an element of $Aff(\bbc )$
and $(a, [{\bf v}, {\bf w}])$ as a rational map according to
(\ref{identification}).  The quotient by this action is the quotient
of $(\bbc^{n}-\{0\})\times (\bbc^m-\{0\})$ by ${\bbc^*}\times
{\bbc^*}$ acting as follows: $(\alpha,a)\times ({\bf v},{\bf
  w})\mapsto ({1\over a\alpha}{\bf v}, a {\bf w})$.  This is
equivalent to the action by multiplication componentwise so that the
quotient is $\brat{1}(\gr )\cong\bbp^{n-1}\times\bbp^{m-1}$ and the
projection $\rat{1}(\gr )\lrar\brat{1}(\gr )$ sends $(a,{\bf v},{\bf
  w})\mapsto ([{\bf v}],[{\bf w}])$. We need to see next how to relate this
  map to the kernel and span. Our first claim is
  that for $f\in\rat{1}(\gr )$,
  $\dim\hbox{ker}(f) = n-1$
  and $\dim\hbox{span}(f)=n+1$. The basing is
$f(\infty)=E_n\subset\bbc^{n+m}$, where by definition $E_n$ is the plane spanned by the
\textit{first} $n$-coordinates vectors. We then have
$\hbox{ker} (f)\subset E_n\subset \hbox{span}(f)$.  The subspace of
$Fl_{(n-1,n+1)}(\bbc^{n+m})$ denoted $Fl_{E_n}$ of all $(A\subset B)$ flags such that
$A\subset E_n\subset B$ is easily identified with $\bbp^{n-1}\times\bbp^{m-1}$ and we have
the map
\begin{equation}\label{kerspan}
\brat{1}(\gr )\lrar \bbp^{n-1}\times\bbp^{m-1}\ \ \ ,\ \ \
f\longmapsto (\hbox{ker}(f)\subset \hbox{span}(f))
\end{equation}
which we claim must now be a homeomorphism. We check the details next.

First we compute the dimension of kernel and span.
To that end we need know that the map $f: z\mapsto {1\over z-a}A$ as
described in (\ref{rational}), with
$A = [v_1,\ldots, v_n]\cdot [w_1,\ldots, w_m]^T$ of rank one,
 can also be viewed as the map sending
$z$ to the $n$-plane $f(z)$ in $\gr$ given as the span of the row vectors in
the matrix
\begin{equation}\label{matrixform}
\begin{pmatrix}
1&0&\cdots&0&{v_1w_1}&\cdots&{v_1w_m}\\
0&1&\cdots&0&{v_2w_1}&\cdots&{v_2w_m}\\
\vdots&&&\vdots&\vdots&&\vdots\\
0&0&\cdots&z-a&{v_nw_1}&\cdots&{v_nw_m}
\end{pmatrix}
\end{equation}
This matrix is well-defined up to the left action by $GL_n(\bbc [z] )$
which corresponds to taking row operations. By dividing up by ${z-a}$
and letting $z\rightarrow +\infty$, we see that indeed $f(\infty
)=E_n\subset\bbc^{n+m}$.  We recall from Remark \ref{description} that
$T(z)={1\over z-a}A$ is the transfer function associated to (\ref{matrixform}).
Out of this matrix representation (\ref{matrixform}),
$\hbox{span}(f)$ becomes by definition the span of the row vectors making
up (\ref{matrixform}) as $z$ varies. This is easily seen to be
the $(n+1)$-dimensional subspace spanned by the $(n+1)$-row vectors
$$\left\{\begin{matrix}1&\cdots&0&0&\cdots&0\\
\vdots&&\vdots&&\vdots\\
0&\cdots&1&0&\cdots&0\\
0&\cdots&0&w_1&\cdots&w_m
\end{matrix}
\right.$$ and is uniquely determined by the homogeneous vector
$[w_1:\cdots :w_m]$.  On the other hand the kernel can be seen to
correspond to
$$\hbox{ker}(f) = \left\{\alpha_1\begin{pmatrix}
1\\
\vdots\\
0\\
\vdots\\
0\end{pmatrix} + \cdots + \alpha_n\begin{pmatrix}
0\\
\vdots\\
1\\
\vdots\\
0\end{pmatrix}\ \ ,\ \ \sum\alpha_iv_i = 0\right\}
$$
and this is an $(n-1)$-dimensional subspace determined by
the homogeneous vector $[v_1:\cdots :v_n]$. The pair
$(\hbox{ker}(f),\hbox{span}(f))$ such that $\hbox{ker} (f)\subset
E_n\subset \hbox{span}(f)$ determines a unique point in
$\bbp^{n-1}\times\bbp^{m-1}$ and this sets up an explicit
homeomorphism with $\brat{1}(\gr )$.

We now turn to the unbased holomorphic maps.
Note that $U(n+m)$ acts on all of $\hol{1}(\gr )$ and
$\bhol{1}(\gr )$ by acting on the target and that the (ker,span)-map
is $U(n+m)$-equivariant. We will construct directly an inverse
map $\phi$ to this (ker, span) map
\begin{equation}\label{phi}
\phi : Fl_{(n-1,n+1)}(\bbc^{n+m})\lrar\bhol{1}(\gr )
\end{equation}
and hence show that $\phi$ is a homeomorphism.
Such a map exists when restricted to the subspace $Fl_{E_n}\cong\bbp^{n-1}\times\bbp^{m-1}$
of flags $(A\subset E_n\subset B)$ as
we've just shown and we write $\phi_r$ such a restriction.
Let $(E_{n-1}\subset E_{n+1})$ be the trivial flag, where the basis vectors for $E_k$ are
taken to be the first $k$-coordinate vectors.
Then $\phi_r (E_{n-1}\subset E_{n+1})$ is the class under
the $Aff(\bbc)$-action of the map $z\longmapsto {1\over
  z}\begin{pmatrix}
  0&\cdots&0\\
  \vdots&&\vdots\\
  1&\cdots&0
\end{pmatrix}
$ corresponding to picking up
$[w_1:\cdots :w_m] = [1:\cdots :0:0]$ and $[v_1:\cdots :v_n]=[0:\cdots :0:1]$.
Equivalently this is the map which in matrix form is
\begin{equation}\label{firstmap}z\longmapsto \hbox{span row vectors}
\begin{pmatrix}
1&\cdots&0&0&0&\cdots&0\\
\vdots&&&&&&\vdots\\
0&\cdots&1&0&0&\cdots &0\\
0&\cdots&0&z&1&\cdots&0\\
\end{pmatrix}
\end{equation}
the righthand side being an $n\times (n+m)$-matrix.

A unitary matrix $g\in U(n+m)$ acts on a flag $(A\subset B)\in Fl_{(n-1,n+1)}(\bbc^{n+m})$
by sending it to the flag $(g(A)\subset g(B))$.
Given any such flag $(A\subset B)$,
there is an element $g$ of $U(n+m)$ sending it to $(E_{n-1}\subset E_{n+1})$.
We can then construct $\phi$ in (\ref{phi}) by setting
\begin{equation}\label{inverse}
\phi (A\subset B) = g^{-1}\phi_r (g(A\subset B)) = g^{-1}\phi_r (E_{n-1}\subset E_{n+1})
\end{equation}
We have to prove this map is well-defined.  Any other element $h\in
U(n+m)$ taking $(A\subset B)$ to the trivial flag satisfies the
property that $hg^{-1}\in U(m-1)\times U(2)\times U(n-1)$ as the
subgroup of $U(n+m)$ consisting of blocks
$\begin{pmatrix}U(n-1)&&\\&U(2)&\\&&U(m-1)\end{pmatrix}$. Different
choices of matrices in $U(n-1)$ and $U(m-1)$ do not affect
$\phi(A\subset B)$ since they stabilize the flag $E_{n-1}\subset E_n\subset E_{n+1}$ and
hence don't affect $\phi_r (E_n\subset
E_{n+1})$. It remains to analyze the effect of $U(2)$ on
(\ref{inverse}). Let $g=\begin{pmatrix}a&b\\c&d\end{pmatrix}\in U(2)$.
This stabilizes the trivial flag but takes $E_{n+1}$ from its standard representation as
span of the first $n+1$-coordinate vectors to the span
$$E'_{n+1} = \hbox{span row vectors}
\begin{pmatrix}
\begin{pmatrix}Id_{n-1}\end{pmatrix}&\\
&a&b&0&\cdots&0\\
&c&d&0&\cdots&0
\end{pmatrix}
$$
Of course $E'_{n+1}=E_{n+1}$ but we distinguish them in notation since a priori the map
$\phi$ can depend on this choice of basis.
According to (\ref{inverse}), $(E_n\subset E'_{n+1})$  must be mapped to the holomorphic
map obtained by acting on each row vector of (\ref{firstmap}) by $g\in
U(2)\subset U(n+m)$. This gives the map
$$
z\longmapsto  \hbox{span row vectors}
\begin{pmatrix}
1&\cdots&0&0&0&\cdots&0\\
\vdots&&&&&&\vdots\\
0&\cdots&1&0&0&\cdots &0\\
0&\cdots&0&az+b&cz+d&\cdots&0\\
\end{pmatrix}
$$
or equivalently the map
$$
z\longmapsto  \hbox{span row vectors}
\begin{pmatrix}
1&\cdots&0&0&0&\cdots&0\\
\vdots&&&&&&\vdots\\
0&\cdots&1&0&0&\cdots &0\\
0&\cdots&0&{az+b\over cz+d}&1&\cdots&0\\
\end{pmatrix}
$$
This is now evidently the composition of the map in (\ref{firstmap})
with the transformation $z\mapsto {az+b\over cz+d}$. This means in summary
that the candidate $\phi (A\subset B)$ we constructed in (\ref{inverse})
is well-defined as an element of $\bhol{1}(\gr )$ and is clearly an inverse
to the (ker,span)-map.

Remains to show this map is continuous but this is
a consequence of the fact that $\phi$ is the induced quotient map in the following
diagram
$$\xymatrix{Fl_{E_n}\times U(n+m)\ar[dr]^{\tilde\phi}\ar[d]^\pi\\
Fl_{(n-1,n+1)}(\bbc^{n+m})\ar[r]^\phi&\bhol{1}(\gr )}
   $$
where $\tilde\phi (f,g) = g\phi_r(f)$ and $\pi (f,g)= gf$.
Theorem \ref{main3} is proved.

\section{Homological Calculations}

This final section uses the sphere bundle description in
Theorem \ref{hol1} to calculate the homology of $\hol{1}(\gr )$ for
various ring coefficients and small values of $n$ and $m$.
We will be working with an explicit description for $H^*(\gr ;\bbz
)$. This cohomology is generated by the Chern classes of the
tautological bundle of $\gr$. Consider the two natural embeddings
$\gr\hookrightarrow BU(n)$ and $\gr\hookrightarrow BU(m)$.  The
cohomology of $\gr$ is generated by the pullbacks of the chern classes
$c_i\in H^*(BU(n))$ and $\bar c_j\in H^*(BU(m))$.  These pullback
classes satisfy the relation
\begin{equation}\label{relation}
(1+c_1+\cdots +c_n)(1+\bar c_1+\cdots +\bar c_m)=1
\end{equation}
This relation is the consequence of the fact that the total chern
classes multiply as in $c(\gamma_m )c(\gamma_n ) =
c(\gamma_m\oplus\gamma_n)$ and that
$\gamma_m\oplus\gamma_n=\gamma_m\oplus Q_m\cong \epsilon_{n+m}$ is
trivial on $\gr$.  The following is standard.

\bpr\label{amodi} There is a graded ring isomorphism
$$H^*(\gr
) = \bbz [c_1,\ldots, c_n, \bar c_1, \ldots, \bar
c_m]/\left(\sum_{i+j=k}c_i\cdot \bar c_j; 1\leq k\leq n+m \right)$$
with $\deg c_i = \deg \bar c_j = 2i$. \epr

\bre\label{relations} The relations $\sum_{i=0}^kc_i\times\bar
c_{k-i}=0$ for $k=1,\ldots, m$ can be solved inductively to express
$\bar c_1, \ldots, \bar c_m$ in terms of $c_1,\ldots c_n$, and the
ring $H^*(\gr )$ is in fact a quotient of $\bbz [c_1,\ldots, c_n]$ by
an explicit ideal $(\rho_1,\ldots, \rho_n)$. That is and for $m\geq
n\geq 2$, $H^*(\gr ) = \bbz [c_1,\ldots, c_n]/(\rho_1,\ldots, \rho_n)$
with $\deg\rho_i = 2m+2i$, $1\leq i\leq n$, and the $\rho_i$ making up
a so-called \textit{regular sequence} (complete formulae in
\cite{chair},\S2).  In particular and when $n=2$, $H^*(Gr(2,m)) = \bbz
[c_1,c_2]/(\rho_1,\ldots, \rho_m)$ with
\begin{equation}\label{formule} \rho_i = \sum_{p_1+2p_2
= m+i}(-1)^{p_1+p_2}{p_1+p_2\choose p_1}c_1^{p_1}c_2^{p_2}
\end{equation}
\ere

We work out explicitly the case of the quadric Grassmann manifold $G(2,2)$.

\bco\label{cohquadric} $H^*(\grass ) = \bbz
[c_1,c_2]/(c^3_1-2c_1c_2, c_1^4 - 2c_2^2)$. \eco

\begin{proof} We ``peel off" the relations $(1+c_1+c_2)(1+\bar c_1+\bar c_2) =
  1$ one by one. In homological degree $2$, $c_1+\bar c_1 = 0$ so that $\bar
  c_1=-c_1$. In degree four, $c_1\bar c_1 + c_2 + \bar c_2 = 0$ so that $\bar
  c_2 = c_2-c_1^2$. In degree six $c_1\bar c_2 + \bar c_1 c_2 = 0$ leading to
  the first relation $c^3_1-2c_1c_2=0$.  The relation $c_1^4=2c_2^2$ follows
  similarly. This is of course consistent with (\ref{formule}).
\end{proof}

More generally we have the following result of Baum needed in
\S\ref{quadric} when computing the cohomology of some flag manifolds.

\bth\label{baum} If $G$ is a compact connected Lie group and $U$ a
closed connected subgroup, both of whose
integral cohomology rings are exterior algebras on
(odd degree) generators, then the sequence
$$H^*(BG)\fract{\rho^*}{\lrar} H^*(BU)
\fract{\sigma^*}{\lrar} H^*(G/U)$$ has the property that the kernel of
$\sigma^*$ is the ideal of $H^*(BU)$ generated by the elements of
positive degree in Image($\rho^*$).
\end{theorem}

\bre $H^*(\gr )$ being a quotient of $H^*(BU(n))$ shows that $H_*(\gr
)$ injects into $H_*(BU(n))$, and thus the embedding
$\gr\hookrightarrow Gr(n,m+1)$ induces a monomorphism in
homology. This fact is no longer true for real Grassmann varieties.
\ere


\subsection{The Quadric Grassmannian}\label{quadric}

In this section we give complete rational calculations for
$\hol{1}(Gr(2,2))$ where $Gr(2,2)$ is the quadric grassmannian.  We do
this by first understanding the Gysin sequence for the bundle $\xi$ derived
from Theorem \ref{hol1}
\begin{equation}\label{bundleg22}
S^3\lrar Hol_1(G(2,2))\lrar Fl(1,2,4)
\end{equation}
This is given by (setting $Fl:=Fl(1,2,4)$ and $G=\grass$ )
\begin{equation}\label{SG}
\cdots\rightarrow H^{i+3}Hol_1G\lrar H^i FL\fract{\chi}{\lrar} H^{i+4}Fl \lrar
H^{i+4}Hol_1G\rightarrow \cdots
\end{equation}
where $\chi$ is the cup-product with the Euler class of $\Fl$; that is
with the top chern class of the bundle $p_1^*(\gamma )\tensor
p_2^*(Q^\vee)$ in Theorem \ref{hol1}.

To begin we need the cohomology structure of $Fl(1,2,4)=
U(4)/U(1)\times U(1)\times U(2)$. Notice that there is a fibration
\begin{equation}\label{fltog}
\bbp^1\lrar Fl(1,2,4)\lrar G(2,2)
\end{equation}
which is trivial in cohomology (additively) since all generators are
concentrated in even degree. This says that $H^*(\Fl )$ is torsion
free with Poincar\'e polynomial
\begin{eqnarray}
{\mathcal P}(\Fl) &=& {\mathcal P}(\bbp^1 )\cdot {\mathcal P}(G(2,2))\nonumber\\
&=&(1+t^2)(1+t^2+2t^4+t^6+t^8)\nonumber\\
&=&1+2t^2+3t^4+3t^6+2t^8+t^{10}\label{ppfl}
\end{eqnarray}
To get the cohomology structure however we have to resort to the result
of Baum explained in Theorem \ref{baum}.

\ble\label{cohg22} There is a ring isomorphism
\begin{eqnarray*}
H^*(\Fl)&=&\bbz [x,y]/((x+y)(x^2+y^2),x^3y+xy^3+x^2y^2)\\
&=&\bbz [x,y]/((x+y)(x^2+y^2),x^4,y^4)
\end{eqnarray*}
where $|x| = 2, |y| = 2$.
\ele

\begin{proof} Consider the short exact sequence
$$H^*(BU(4))\fract{\rho^*}{\ra 3} H^*(BU(1)\times BU(1)\times
BU(2))\fract{\sigma^*}{\ra 3} H^*(\Fl)$$ with kernel of $\sigma^*$ the
ideal of positive degree elements in $\rho^*$.  Here $H^*(BU(1)\times
BU(1)\times BU(2)) = \bbz [x,y,z_1,z_2]$ so that
$$H^*(\Fl ) = \bbz [x,y,z_1,z_2]/I$$
where $I$ is generated by $(1+x)(1+y)(1+z_1+z_2)$; i.e with the ideal
of relations generated by
$$1+(x+y+z_1)+(z_2+xy+xz_1+yz_1)+(yz_2+xz_2+xyz_1)+xyz_2=1$$
From this we can deduce that $z_1=-x-y$, $z_2=-xz_1-xy-yz_1=x^2+y^2+xy$
and that
$(x+y)(x^2+y^2)=0$ et $x^3y+xy^3+x^2y^2=0$ as claimed.
The second isomorphism follows by identifying generators in each
dimension as in (\ref{ppfl}).
\end{proof}

\bre \label{cohoeffect} As we have indicated there are two maps
$$p_1 : \Fl\lrar \bbp^3\ \ \ \hbox{and} \ \ \ p_2 : \Fl\lrar Gr(2,2)$$
and so we determine their effect on cohomology. We write $H^*(Gr(2,2))
= \bbz [c_1,c_2]/I$ and $H^*(\bbp^3) = \bbz [u]/(u^4)$. The first map
is modeled after
$$U(4)/U(1)\times U(1)\times U(2)\lrar U(4)/U(1)\times
U(3)$$
which is induced from the multiplication $U(1)\times U(2)\lrar U(3)$.
There is a diagram
$$\xymatrix{ \Fl\ar[r]\ar[d]&BU(1)\times BU(1)\times
  BU(2)\ar[d]^\pi\\ \bbp^3\ar[r]&BU(1) }$$ where $\pi$ is projection
onto the first component.  Since the generator $u$ is induced from
that single copy of $BU(1)$, $f^*(u) = x$. Similarly $p_2$ is
determined from multiplication $H^*(BU(2))\lrar H^*(BU(1)\times
BU(1))$ which sends $1+c_1+c_2$ to $(1+x)(1+y)$ so that
$$p_2^*(c_1) = x+y\ \ \ ,\ \ \ p_2^*(c_2) = xy$$
\ere

We are now in a position to read off the Euler class of the bundle
(\ref{bundleg22}), or equivalently the top chern class $c_2$ of
$\xi := p_1^*(\gamma )\tensor p_2^*(Q^\vee)$.

\ble\label{euler} $c_2(\xi)=3x^2+y^2+2xy\in H^4(Fl(1,2,4))=\bbz \{x^2,xy,y^2\}$.
\ele

\begin{proof} By the formula for the total chern class of a tensor
  product of two bundles $L,E$ with $L$ a line bundle, one gets the
  general relation $c_2(E\otimes L)= c_2(E)+c_1(E)\times
  c_1(L)+c_1^2(L)$ (see \cite{fulton}, p.55). Applying this to our
  situation with $L=p_1^*(\gamma )$ and $E=p_2^*(Q^\vee )$, we get
$$c_2(\xi ) =(x^2+y^2+xy)+x(x+y)+x^2 = 3x^2+y^2+2xy$$
as claimed.
\end{proof}

An analysis of the Gysin sequence (\ref{SG}) with integral coefficients yields
the following complete calculation stated in Theorem \ref{grasscase}.

\bth\label{homcalc} $\tilde H^*(\hol{1}(\grass );\bbz )$ is given by
$$
\begin{tabular}{|c|c|c|c|c|c|c|c|c|c|c|c|c|c|c|c|c|c|c|c|}
  \hline
  $i$ & 2
  & 4 & 6 & 7 & 8 & 9 & 11 & 13 \\
  $H^i$& $\bbz\oplus\bbz$ & $\bbz\oplus\bbz$ & $\bbz_4\oplus\bbz$ & $\bbz$ & $\bbz_4$ & $\bbz\oplus\bbz$ & $\bbz\oplus\bbz$ & $\bbz$ \\ \hline
\end{tabular}$$
and $0$ in all other degrees.
\end{theorem}

\begin{proof}
In the Gysin sequence (\ref{SG}), the map $\chi$ is multiplication by
$3x^2+y^2+2xy$. Since the cohomology of $\Fl$ is concentrated in
even dimension, the exact sequence splits into short exact sequences
$$
0\rightarrow H^{2i+1}Hol_1G(2,2)\rightarrow H^{2i-2} \Fl\fract{\chi}{\lrar}
  H^{2i+2}\Fl\rightarrow H^{2i+2}Hol_1G(2,2)\rightarrow 0 $$ from which we deduce
  that
\begin{eqnarray*}
H^{odd}Hol_1(G(2,2)) &=& \hbox{ker}( H^{odd-3}\Fl\lrar H^{odd+1}\Fl)\\
H^{even}Hol_1(G(2,2)) &=& coker( H^{even-4}\Fl\lrar H^{even}\Fl)
\end{eqnarray*}
Now we list generators for $H^*(\Fl )$
\begin{eqnarray}
  H^2 &=& \bbz \{x,y\}\nonumber\\
  H^4 &=& \bbz \{x^2, xy, y^2\}\nonumber\\
  H^6 &=& \bbz \{x^3,x^2y, y^2x\}\ \ \ y^3=-x^3-x^2y-xy^2 \label{generators}\\
  H^8 &=& \bbz \{x^3y,y^3x\}\ \ \ \ \ \ \ x^4=y^4=0 \ \ ,\ \ x^2y^2=-xy^3-x^3y\nonumber\\
  H^{10} &=& \bbz \{x^3y^2\}\ \ \ \ \ \ \ \ \ \ \ \ x^3y^2=-x^2y^3\nonumber
\end{eqnarray}
where here $H^* = H^*(\Fl )$ and we've written the relations on the right for
each degree, and then analyze multiplication by the euler class on each.
It is clear that $H^1\Hl=0$ and $H^2\Hl=\ H^2\Fl=\bbz\{x,y\}$.
To compute $H^5$ and $H^6$ for example we consider the portion of the sequence
$$0\rightarrow H^5\Hl\rightarrow H^2\Fl\fract{\chi}{\lrar} H^6\Fl\rightarrow  H^6\Hl\rightarrow 0$$
with $H^2\Fl = \bbz\{x,y\}$ and $H^6\Fl = \bbz\{x^3,x^2y,xy^2\}$.  In the
bases $(x,y)$ and $(x^3,x^2y,xy^2)$, the matrix of $\chi$ is $\left(
  \begin{array}{cc}
    3 & -1 \\
    2 & 2 \\
    1 & 1 \\
  \end{array}
\right)$ with Smith associated form
$\left(
  \begin{array}{cc}
    1 & 0 \\
    0 & 4 \\
    0 & 0 \\
  \end{array}
\right)$. Thus the cokernel is $\bbz\oplus \bbz_4$ and
$H^5\Hl=0$.
Similarly
$$0 \rightarrow H^7\Hl\rightarrow H^4\Fl\fract{\chi}{\lrar} H^8\Fl\rightarrow
H^8\Fl\rightarrow 0$$ with $H^4\Fl=\bbz^3\{x^2,xy,y^2\}$,
$H^8\Fl=\bbz^2\{x^3y,xy^3\}$ and $x^2y^2=-x^3y-xy^3$.
Keeping in mind the relation $x^2y^2=-x^3y-xy^3$ we see that
\begin{eqnarray*}
\chi (x^2) &=& 3x^4+x^2y^2+2x^3y=-x^3y-xy^3+2x^3y=x^3y-xy^3\\
\chi(xy)&=&3x^3y+xy^3+2x^2y^2=3x^3y+xy^3+2(-x^3y-xy^3)\\ &&=x^3y-xy^3\\
\chi(y^2)&=&3x^2y^2+y^4+2xy^3=3(-x^3y-xy^3)+2xy^3 =-3x^3y-xy^3
\end{eqnarray*}
The matrix of
$\chi$ is $\left(
  \begin{array}{ccc}
    1  & 1  &-3 \\
    -1 & -1 &-1  \\
  \end{array}
\right)$ with associated Smith form
 $\left(
  \begin{array}{ccc}
    1  & 0  & 0  \\
    0  & 4  & 0 \\
  \end{array}
\right)$
and hence  $H^7\Hl=\bbz$  et  $H^8\Hl=\bbz_4\oplus\bbz$. The rest
follows similarly.
\end{proof}

\bco\label{CohG22}  The rational
Poincar\'e series for the holomorphic mapping space is
$${\mathcal P}(\Hl) =1+2t^2+2t^4+t^6+t^7+2t^9+2t^{11}+t^{13}$$
\eco

\bre From (\ref{bundleg22}) (or corollary \ref{unitsphere}), we know that
$\hol{1}(\grass )$ is of the homotopy type of a closed manifold of real
dimension $13$. This explains the Poincar\'e duality for the coefficients
appearing in Corollary \ref{CohG22}.  \ere


\subsection{The general case $(1,2,n+2)$-flags} We try to extend the results
on the quadric grassmannian to maps into $Gr(2,n)$.  We generalize lemma
\ref{cohg22} first

\ble\label{cohflag} We have the algebra isomorphism
$$
H^*Fl(1,2,n+2)=\bbz[x,y]/\left(\sum_{i+j=n+1}
x^iy^j,\ x^{n+2}=0=y^{n+2}\right) $$
\ele

\begin{proof} As in the proof of Lemma \ref{cohg22},
$H^*Fl(1,2,n+2)=\bbz[x,y]/I$ where $I$ is the ideal generated by all
homogeneous relations in $x$ and $y$ derived from
$$(1+x)(1+y)(1+z_1+z_2+\cdots+z_n)=1$$
These we get recursively
$$(xy)\ z_{k-2}\ +\ (x+y)\ z_k\ +\ z_k=0\ ,\ \ \ \forall
  k>1$$
We can write this down in matrix form
$$\left(
\begin{array}{ccc}
z_k \\ z_{k-1}
\end{array}
\right)=
\left(
\begin{array}{ccc}
-(x+y) & -xy  \\
1 & 0
\end{array}
\right)
\left(
\begin{array}{ccc}
z_{k-1}   \\
z_{k-2}
\end{array}
\right)=\left(
\begin{array}{ccc}
-(x+y) & -xy  \\
1 & 0
\end{array}
\right)^{k-1}.
\left(
\begin{array}{ccc}
-(x+y)   \\
1
\end{array}
\right)
$$
After diagonalizing we obtain the power matrix
$$\left(
\begin{array}{ccc}
-(x+y) & -xy  \\
1 & 0
\end{array}
\right)^{k}=\frac{1}{y-x} \left(
\begin{array}{ccc}
(-x)^{k+1}-(-y)^{k+1} &  y(-x)^{k+1}-x(-y)^{k+1} \\
(-x)^{k}-(-y)^{k} & y(-x)^{k}-x(-y)^{k}
\end{array}
\right)$$
from which we deduce all relations
$$z_n=(-1)^{n}\displaystyle{\sum_{i+j=n}}x^iy^j$$
$$z_{n+1}\ =\ 0 \ =\ -(xy)\ z_{n-1}\ -\ (x+y)\ z_{n}\ \Longrightarrow\
\displaystyle{\sum_{i+j=n+1}}x^iy^j=0$$
$$z_{n+2}\ =\ 0 \ =\ -(xy)\ z_{n}\ -\ (x+y)\ z_{n+1}\ \Longrightarrow\
xy\displaystyle{\sum_{i+j=n}}x^iy^j=0$$
so that
$$H^*Fl(1,2,n+2)=\bbz[x,y]/\left(\sum_{i+j=n+1}
x^iy^j,\ xy\displaystyle{\sum_{i+j=n}}x^iy^j\right)$$
Now we need see that both ideals
 $$\left(\sum_{i+j=n+1}
x^iy^j,\ xy\displaystyle{\sum_{i+j=n}}x^iy^j\right)\ \leftrightarrow
\ \left(\sum_{i+j=n+1} x^iy^j,\ x^{n+2}, y^{n+2}\right) $$ are equivalent
but this is immediate.
\end{proof}

Next we compute the euler class $\mathbf{e}(\xi_n)$ of the sphere fibration
(Theorem \ref{hol1})
$$\xi_n\ :\ S^{2n-1}\lrar Hol_1(G(2,n))\lrar Fl(1,2,n+2)$$

\bpr\label{eulerg2n}
$\displaystyle\mathbf{e}(\xi_n )=\frac{\partial}{\partial x}\displaystyle{\sum_{i+j=n+1}}x^iy^j$.
\epr

\begin{proof}
This euler class corresponds to the top chern class
$c_n(p_1^*H\tensor p_2^*Q^\vee)$
of $\zeta (\xi_n)$ and this can be computed by the formula
$${\bf e}(\xi_n ) = \displaystyle{\sum_{i=0}^n} c_1^i(p_1^*H)\cdot
c_{n-i}(p_2^*Q^\vee)$$
We have to determine the chern classes of the pullback of $Q^\vee$.
A quick inspection of the diagram
$$\xymatrix{ Fl(1,2,n+2)
  \ar[d]\ar[r]^{p_2}&Gr(2,n)\ar[d]\\ BU(1)\times
  BU(1) \times BU(n)\ar[r]& BU(2) \times BU(n) }$$
  shows that the map in cohomology
$H^*(Fl(1,2,n+2))\longleftarrow
H^*(Gr(2,n))=\bbz[c_1,c_2]/I$, with $I$ the ideal generated by the
$\rho_1,\ldots, \rho_n$ given in (\ref{formule}), sends
\begin{eqnarray*}
c_1 & \longmapsto & x+y\\
c_2 & \longmapsto & xy\\
\overline{c}_k = c_k(Q )
& \longmapsto & z_n=(-1)^{k}\displaystyle{\sum_{i+j=k}}x^iy^j
\end{eqnarray*}
using the fact that the multiplication map $BU(1)\times BU(1)\lrar BU(2)$
in cohomology sends $c_1\mapsto x+y$ and $c_2\mapsto xy$.
But $c_{k}(Q^\vee )=\ (-1)^k \overline{c}_k(Q)$ so that
$$p_2^*(c_k(Q^\vee )) = (-1)^{2k}
\displaystyle{\sum_{i+j=k}}x^iy^j=\ \displaystyle{\sum_{i+j=k}}x^iy^j
\in H^*(Fl(1,2,n+2))$$
To show that $\mathbf{e}\ (\xi_n)=
\frac{\partial}{\partial x}\displaystyle{\sum_{i+j=n+1}}x^iy^j$, we proceed by induction.
The case $n=2$ has been verified in lemma \ref{euler}
where we checked that $\mathbf{e}(\xi_2 )=3x^2+2xy+y^2$.
The proof now proceeds inductively
\begin{eqnarray*}
\mathbf{e}\ (\xi_{n+1}) =\  c_{n+1}(Q^\vee)+ \mathbf{e}\ (\xi_n)
&=& \displaystyle{\sum_{i+j=n+1}}x^iy^j+ \frac{\partial}{\partial x}\displaystyle{\sum_{i+j=n+1}}x
^iy^j\\
&=& \displaystyle{\sum_{i+j=n+1}}x^iy^j+ \displaystyle{\sum_{i+j=n+1}}i\ x^iy^j\\
&=& \frac{\partial}{\partial x}\displaystyle{\sum_{i+j=n+2}}x^iy^j
\end{eqnarray*}
which is what we wanted to prove.
\end{proof}

As an application we can give the explicit ranks for
$H^*(Fl(1,2,n+2);\bbq )$ as well as the euler class of the sphere
bundle $\xi$. It is then possible by use of the Gysin sequence again
to obtain complete calculations for
$H^*(\hol{1}(Gr(2,n))$ with all field coefficients. This was done for
the quadric grassmannian $\grass$ in Corollary \ref{CohG22}.

\bco\label{CohG23}. The rational Betti numbers for $\hol{1}(Gr(2,3))$
are listed below
$$
\begin{tabular}{|c|c|c|c|c|c|c|c|c|c|c|c|c|c|c|c|c|c|c|c|}
  \hline
  $i$ & 1
  & 2 & 3 & 4 & 5 & 6 & 7 & 8 & 9 & 10 & 11 & 12 & 13 & 14 & 15 & 16 &
  17 & 18 & 19\\ $b_i $& 0 & 2 & 0 & 3 & 0 & 3 & 0 & 2 & 1 & 1 & 2 & 0
  & 3 & 0 & 3 & 0 & 2 & 0 & 1 \\ \hline
\end{tabular}$$
\eco

As expected these betti numbers satisfy Poincar\'e duality.

\addcontentsline{toc}{section}{Bibliography}
\bibliography{biblio}
\bibliographystyle{plain[8pt]}

\end{document}